\documentclass [a4,10pt,twoside]{article}
\usepackage{amsmath, amssymb, amsthm, amsfonts, amsxtra, latexsym, amscd,
pb-diagram,  graphics,rotating,xy}

\theoremstyle{plain}
\newtheorem{dl}{Theorem}[section]

\newtheorem{dn}[dl]{Definition}

\newtheorem{thm}{\bf Theorem}
 
\newtheorem{pro}[thm]{\bf Proposition} 

\newcommand{\A}{\mathcal{A}}
\newcommand{\La}{\mathcal{L}}
\newcommand{\R}{\mathcal{R}}
\newcommand{\tx}{\otimes}
\newcommand{\ts}{\oplus}
\begin{document}
\title{ON THE AXIOMATICS OF ANN-CATEGORIES}
\author{Nguyen Tien Quang, D. D. Hanh and N. T. Thuy}
\pagestyle{myheadings} 
\markboth{On the axiomatics of Ann-categories}{N. T. Quang, D. D. Hanh and N. T. Thuy}
\maketitle
\setcounter{tocdepth}{1}
\begin{abstract}
In this paper, we have studied the axiomatics of {\it Ann-categories} and {\it categorical rings.} These are the categories with distributivity constraints whose axiomatics are similar with those of ring structures. The main result we have achieved is proving the independence of the axiomatics of Ann-category definition. And then we have proved that after adding an axiom into the definition of categorical rings, we obtain the new axiomatics which is equivalent to the one of Ann-categories. 
\end{abstract}

\section{Introduction}
The definition of {\it Ann-categories} was presented by Nguyen Tien Quang in 1987 [4], which was regarded as a categoricalization of {\it ring} structures. Each Ann-category is charactered by 3 invariants: the ring $R=\Pi_0(\A)$ of classes of invertible objects of $\A,$ the $R$-bimodule $\Pi_1(\A)=\text{Aut}(0)$ and the element $f\in H^3(R, M)$ in the third cohomology group of the ring $R$ with coefficients in the $R$-bimodule $M$ due to [5]. Recently, we have proved that this cohomology coincides with the one due to Maclane[3]. The class of regular Ann-categories (whose commutativity constraints satisfy the condition $c_{X,X}=id$) is classificated by the cohomology group $H^3_M(R,M)$ of the $Z$-algebra due to Shukla[8].

In [1], M.Jibladze and T. Pirashvili presented the definition of {\it categorical rings} as a slightly modified version of the definition of Ann-categories and classificated them with Maclane cohomology for rings.

In this paper, we have made some comments on these two definitions. First, we have proved that in the axiomatics of Ann-categories, the compatibility of the functors $(L^A,\widetilde{L}^A)$, $(R^A,\widetilde{R}^A)$ with the commutativity constraint $c$ is dependent. So we have proved that each Ann-category is a categorical ring due to [1]. We have seen that, in order to prove the converse, we must add an axiom into the definition of categorical rings, that is the compatibility of the functors $(L^A,\widetilde{L}^A),(R^A,\widetilde{R}^A)$ with the unitivity constraint. In [6], [7], thanks to this axiom, we may construct the associative bimodule structure. There is a problem here: Is the new added axiom independent of the others in the definition of categorical rings due to [1]?

For convenience, in this paper we denote by $AB$ the tensor product of the two objects $A$ and $B,$ instead of $A\tx B.$

\section{The axiomatics of Ann-categories}
\begin{dn} An Ann-category consists of:\\
i) A groupoid $\A$ together with two bifunctors $\ts,\tx:\A\times\A\longrightarrow \A.$\\
ii) A fixed object $0\in \A$ together with naturality constraints $a^+,c,g,d$ such that $(\A,\ts,a^+,c,(0,g,d))$ is a Pic-category.\\
iii) A fixed object $1\in\A$ together with naturality constraints $a,l,r$ such that $(\A,\tx,a,(1,l,r))$ is a monoidal $A$-category.\\
iv) Natural isomorphisms $\La,\R$
\[\La_{A,X,Y}:A\tx(X\ts Y)\longrightarrow (A\tx X)\ts(A\tx Y)\]
\[\R_{X,Y,A}:(X\ts Y)\tx A\longrightarrow(X\tx A)\ts(Y\tx A)\] 
such that the following conditions are satisfied:\\
(Ann-1) For each $A\in \A,$ the pairs $(L^A,\breve{L^A}),(R^A,\breve{R^A})$ determined by relations:\\
\[
\begin{aligned}
&L^A & = &A\tx- \;\;\;\; &R^A&=&-\tx A\\
&\breve{L^A}_{X,Y}& = &\La_{A, X, Y}\;\;\;\; &\breve{R^A}_{X,
Y}&=&\mathcal{R}_{X, Y, A}
\end{aligned}
\]
are $\ts$-functors which are compatible with $a^+$ and $c.$\\
(Ann-2) For all $ A,B,X,Y\in \A,$ the following diagrams:

\[ 
\begin{diagram}
\node{(AB)(X\ts Y)}\arrow{s,l}{\breve{L}^{AB}} \node{ A(B(X\ts
Y))}\arrow{w,t}{\quad a_{A, B, X\ts Y}\quad}\arrow{e,t}{\quad
id_A\tx\breve{L}^B\quad}\node{ A(BX\ts
BY)}\arrow{s,r}{\breve{L}^A}\\
\node{(AB)X\ts(AB)Y}\node[2]{A(BX)\ts A(BY)}\arrow[2]{w,t}
{\quad\quad a_{A, B, X}\ts a_{A, B, Y} \quad\quad}
\end{diagram}\tag{1.1} 
\]

\[ \begin{diagram}
\node{(X\ts Y)(BA)}\arrow{s,l}{\breve{R}^{BA}}\arrow{e,t}{\quad
a_{ X\ts Y, B, A}\quad}\node{((X\ts Y)B)A}\arrow{e,t}{\quad
\breve{R}^B\tx id_A\quad}\node{ (XB\ts
YB)A}\arrow{s,r}{\breve{R}^A}\\
\node{X(BA)\ts Y(BA)}\arrow[2]{e,t} {\quad\quad a_{X, B, A}\ts
a_{Y, B, A} \quad\quad}\node[2]{(XB)A\ts (YB)A}
\end{diagram}\tag{1.1'} \]
\[ \begin{diagram}
\node{(A(X\ts Y))B}\arrow{s,l}{\breve{L}^{A}\tx id_B} \node{
A((X\ts Y)B)}\arrow{w,t}{\quad a_{A, X\ts Y,
B}\quad}\arrow{e,t}{\quad id_A\tx\breve{R}^B\quad}\node{
A(XB\ts
YB)}\arrow{s,r}{\breve{L}^A}\\
\node{(AX\ts AY)B}\arrow{e,t}{\quad \breve{R}^B
\quad}\node{(AX)B\ts (AY)B}\node{A(XB)\ts A(YB)}\arrow{w,t}
{\quad a\ts a\quad}
\end{diagram}\tag{1.2} \]
\[\begin{diagram}
\node{(A\ts B)X\ts(A\ts
B)Y}\arrow{s,r}{\breve{R}^{X}\ts\breve{R}^Y} \node{(A\ts
B)(X\ts Y)}\arrow{w,t}{\breve{L}^{A\ts B}}\arrow{e,t}{
\breve{R}^{X\ts
Y}}\node{A(X\ts Y)\ts B(X\ts Y)}\arrow{s,l}{\breve{L}^A\ts\breve{L}^B}\\
\node{(AX\ts BX)\ts(AY\ts BY)}\arrow[2]{e,t} {\quad\quad v
\quad\quad}\node[2]{(AX\ts AY)\ts(BX\ts BY)}
\end{diagram}\tag{1.3} \]

\noindent commute, where $v=v_{U,V,Z,T}:(U\ts V)\ts(Z\ts T)\longrightarrow(U\ts Z)\ts (V\ts T)$ is the unique functor built from $a^+,c,id$ in the monoidal symmetric category $(\A,\ts).$\\
(Ann-3) For the unity object $1\in \A$ of the operation $\ts,$ the following diagrams:

\[
\begin{diagram}
\node{1(X\oplus Y)} \arrow[2]{e,t}{\breve{L}^1}
\arrow{se,b}{l_{X\oplus Y}}
\node[2]{1X\oplus 1Y} \arrow{sw,b}{l_X\oplus l_Y} \\
\node[2]{X\oplus Y}
\end{diagram}\tag{1.4}
\]

\[
\begin{diagram}
\node{(X\oplus Y)1} \arrow[2]{e,t}{\breve{R}^1}
\arrow{se,b}{r_{X\oplus Y}}
\node[2]{X1\oplus Y1} \arrow{sw,b}{r_X\oplus r_Y} \\
\node[2]{X\oplus Y}
\end{diagram}\tag{1.4'}
\]

commute.
\end{dn}

{\bf Remark.} The commutative diagrams (1.1), (1.1') and (1.2), respectively, mean that:\\
\[\begin{aligned} (a_{A, B, -})\;:\;& L^A.L^B &\longrightarrow \;& L^{AB}\\
(a_{-,A,B})\;:\;&R^{AB}&\longrightarrow \;&R^A.R^B\\
(a_{A, - ,B})\;:\;&L^A.R^B&\longrightarrow \;&R^B.L^A
\end{aligned}\]
are $\ts$-functors.\\
The diagram (1.3) shows that the family $(\breve{L}^Z_{X,Y})_Z=(\mathcal{L}_{-,X,Y})$ is a $\ts$-functor between the $\ts$-functors $Z\mapsto Z(X\ts Y)$ and $Z\mapsto ZX\ts ZY$, and the family $(\breve{R}^C_{A,B})_C=(\mathcal{R}_{A, B,-})$
is a $\ts$- functor between the functors $C\mapsto (A\ts B)C$ and $C\mapsto AC\ts BC.$\\
The diagram (1.4) (resp. (1.4')) shows that $l$ (resp. $r$) is a $\ts$-functor from $L^1$ (resp. $R^1$) to the unitivity functor of the $\ts$-category $\A$.
\begin{pro}
In the Ann-category $\A$ there exist uniquely isomorphisms:
\[\widehat{L}^A\;:\;A\otimes 0\longrightarrow 0\;,\quad
\widehat{R}^A\;:\;0\otimes A\longrightarrow 0\] such that the following diagrams: 
\[\begin{CD} AX@< L^A(g)<< A(0\ts X)\\
@AgAA   @VV\breve{L}^A V\\
0\ts AX @<\widehat{L}^A\ts id<<A0\ts AX\end{CD}\quad\quad (1.5)\quad\quad\quad
 \begin{CD}AX@< L^A(d)<< A(X\ts 0)\\
@AdAA   @VV\breve{L}^A V\\
AX\ts 0 @<id\ts \widehat{L}^A<<AX\ts A0\end{CD}\tag{1.5'}
\]
\[\begin{CD} AX@< R^A(g)<< (0\ts X)A\\
@AgAA   @VV\breve{R}^A V\\
0\ts AX @<\widehat{R}^A\ts id<<0A\ts XA\end{CD}\quad\quad (1.6)\quad\quad\quad
\begin{CD} AX@< R^A(d)<< (X\ts 0)A\\
@AdAA   @VV\breve{R}^A V\\
AX\ts 0 @<id\ts \widehat{R}^A<<XA\ts 0A\end{CD}\tag{1.6'}
\]
commute.\\
Meaningly, $L^A$ and $R^A$ are functors which are compatible with the unitivity constraint of the operation $\ts$.
\end{pro}
\begin{proof} Since the pair $(L^A, \breve{L}^A)$ is a $\ts$-functor which is compatible with the associativity constraint $a^+$ of the Picard category $(\A,\ts),$ it is also compatible with the unitivity constraint $(0,g,d)$ thanks to Prop.0.4.4 [6]. That means there exists uniquely the isomorphism $\widehat{L}^A$ satisfying the diagrams (1.5) and (1.5').\\
The proof for $\widehat{R}^A$ is completely similar.
\end{proof}

\section{ A remark on the axiomatics of Ann-categories}

The commutativity constraint $c$ plays a quite special role in the study of categories with tensor product. For example, in 1981, Kasangian Stefano and Rossi Fabio [17] presented the problem of the research on the relationship between some conditions for commutativity constraint in symmetric monoidal categories.\\
 
We now consider the axiomatics of Ann-categories in another view. In the definition of a ring as well as a module, the axiom about the commutation of the addition can be omitted: It can be implied from the other axioms. Consider the axiomatics of an Ann-category, we can determine the commutativity constraint $c$ based on the constraints $\mathcal{L}, \mathcal{R}$  and $a^+$ thanks to the commutative diagram (1.3). It leads us to consider the independence or dependence of the axioms related to the commutativity constraint $c.$ That is the compatibility of $c,$ the compatibility of $c$ with $a^+$ and the compatibility of the functors $L^A=A\otimes -, R^A=-\otimes A$ with $c$. In this section, we will prove the independence of the last requirement.

\begin{pro} In the Ann-category $\A,$ the compatibility of the functors $(L^A, \breve{L}^A)$, $ (R^A, \breve{R}^A)$ with the commutativity constraint can be deduced from the other axioms refered in Definition 1, without the compatibility of $c$ and the compatibility of $c$ with $a^+$.

\end{pro}
\begin{proof} First, we prove that the diagram:

\[\begin{CD}
X(A1\oplus B1)@>\mathcal{L}>>X(A1)\oplus X(B1)\\
@Vid_X\otimes c VV @VVcV\\
X(B1\oplus A1)@>\mathcal{L}>>X(B1)\oplus X(A1)
\end{CD}\tag{2.1}\]
commute in step by step:

1. Consider the diagram (2.2), we can see that:

The regions (I), (IV) commute thanks to the naturality of the isomorphism  $\mathcal{L}$,
the regions (II), (VIII), (IX) commute thanks to the axiom (1.1);
the regions (III), (VII) commute thanks to the axiom (1.2);
the regions (V), (X) commute thanks to the axiom (1.3);
the region (VI) commutes thanks to the naturality of the isomorphism $v$.
Therefore, the outside region commutes.

2. Consider the diagram (2.3) in which the region (II) is exactly the outside region of the diagram (2.2) whose commutation was proved right above. 
The religions (I) and (III) commute thanks to the coherence for $\otimes$-functor $(L^X, \breve{L}^X)$; 
 tö $(L^X, \breve{L}^X)$; 
the regions (IV) and (V) commute thanks to the axiom (1.3) and the definition of the isomorphism $v$.
Therefore, the outside region commutes.

3. Now we consider the diagram (2.4) whose outside region is the one of the diagram (2.3).
In this diagram, the religions (I), (II) commute thanks to the naturality of $\mathcal{L}$, so the region (III) commutes. 
Hence, from the regular property of the object $X(A1)$ and $X(B1)$ for the addition $\oplus$, we can deduce that the diagram (2.1) commutes.\\

4. Finally, we prove that the diagram 

\[\begin{CD}
X(A\oplus B)@>\mathcal{L}>>XA\oplus XB\\
@V id\otimes c VV   @V c VV\\
X(B\oplus A)@>\mathcal{L}>>XB\oplus XA
\end{CD}\]
commutes by embedding it into the diagram:

\[\begin{diagram} \node{X(A1\oplus
B1)}\arrow{se,t}{id\otimes(r\oplus r)\quad\quad\quad
\text{(I)}}\arrow[3]{s,r}{id\otimes c\quad
\text{(IV)}}\arrow[3]{e,t}{\mathcal{L}}\node[3]{X(A1)\oplus
X(B1)}\arrow{sw,t}{id\otimes r\oplus id\otimes r}\arrow[3]{s,r}{c}\\
\node[2]{X(A\oplus B)}\arrow{s,r}{id\otimes c\quad \text{(II)}}\arrow{e,t}{\mathcal{L}}\node{XA\oplus XB}\arrow{s,r}{c\quad\quad \text{(V)}}\\
\node[2]{X(B\oplus A)}\arrow{e,t}{\mathcal{L}}\node{XB\oplus XA}\\
\node{X(B1\oplus A1)}\arrow{ne,b}{id\otimes(r\oplus r)\quad\quad
\text{(III)}}\arrow[3]{e,t}{\mathcal{L}}\node[3]{X(B1)\oplus
X(A1)}\arrow{nw,b}{id\otimes r\oplus id\otimes r}
\end{diagram}\]
In this diagram, the outside region is exactly the one of the commutative diagram (2.1), the regions (I) and (III) commute thanks to the naturality of the isomorphism $\mathcal{L}$;
the regions (IV) and (V) commute thanks to the naturality of the isomorphism $c$. So the region (II) commutes.\\
\indent Because of the symmetry, we can deduce the  compatibility of the functor $(R^A,\breve{R}^A)$ with $c$. This completes the proof.

\begin{center}
\setlength{\unitlength}{0.52cm}
\begin{picture}(24,28)
\scriptsize{
\put(7,0){$X(A1\oplus A1)\oplus X(B1\oplus B1)$}
\put(13,0.5){\vector(3,1){3.5}}\put(16,1){$\mathcal{L}\oplus \mathcal{L}$}

\put(0,2){$X((A1\oplus A1)\oplus (B1\oplus B1))$}
\put(4,1.8){\vector(3,-1){3.5}}
\put(5,1){$\mathcal{L}$}
\put(14.8,2){$(X(A1)\oplus X(A1))\oplus(X(B1)\oplus X(B1))$}
\put(16.5,2.5){\vector(-1,1){1.3}}
\put(11.5,3){$(a\oplus a)\oplus (a\oplus a)$}

\put(11,4){$((XA)1\oplus (XA)1)\oplus ((XB)1\oplus (XB)1)$}

\put(11,6){$(XA)(1\oplus 1)\oplus (XB)(1\oplus 1)$}
\put(13.5,5.8){\vector(1,-1){1.3}}
\put(12,5){$\mathcal{L}\oplus \mathcal {L}$}

\put(4,8){$(X(A(1\oplus 1))\oplus(X(B(1\oplus 1))$}
\put(10.2,7.8){\vector(2,-1){2.5}}\put(10,7){$a\oplus a$}
\put(8.5,7.8){\vector(0,-1){7.2}}\put(3.8,4){$(id\otimes \mathcal{L})\oplus(id\otimes \mathcal{L})$}

\put(9,10){$(XA\oplus XB)(1\oplus 1)$}
\put(12.2,9.8){\vector(1,-2){1.6}}\put(13.5,8){$\mathcal{R}$}
\put(9.5,10.7){\vector(0,1){5.2}}\put(9.7,12){$\mathcal{L}$}
\put(0,10){$X(A(1\oplus 1)\oplus B(1\oplus 1))$}
\put(1,9.8){\vector(0,-1){7.3}}\put(1.2,7){$id \otimes (\mathcal{L}\oplus \mathcal {L})$}
\put(4.2,9.8){\vector(3,-1){3.5}}\put(5,9){$\mathcal{L}$}

\put(4,12){$(X(A\oplus B))(1\oplus 1))$}
\put(6.6,11.8){\vector(2,-1){2.5}}\put(5.8,11){$\mathcal{L}\otimes id$}
\put(5,12.7){\vector(0,1){5}}\put(4.3,15.5){$\mathcal{L}$}

\put(10,14){$((XA)1\oplus (XB)1)\oplus ((XA)1\oplus (XB)1)$}
\put(18,13.8){\vector(0,-1){9}}\put(18.3,11){$v$}
\put(0,14){$X((A\oplus B)(1\oplus 1))$}
\put(2,13.8){\vector(2,-1){2.5}}\put(2.7,13){$a$}
\put(1,13.8){\vector(0,-1){3.2}}\put(1.3,11.5){$id\otimes \mathcal{R}$}
\put(1,14.7){\vector(0,1){11}}\put(1.2,21){$id\otimes \mathcal{L}$}

\put(6,16){$(XA\oplus XB)1\oplus (XA\oplus XB)1$}
\put(10.2,15.8){\vector(3,-2){2.2}}\put(12,15){$\mathcal{R}\oplus \mathcal{R}$}

\put(12.4,18){$(X(A1)\oplus X(B1))\oplus (X(A1)\oplus X(B1))$}
\put(18,17.8){\vector(-1,-1){3.2}}\put(13,17){$(a\oplus a)\oplus (a\oplus a)$}
\put(20.5,17.8){\vector(0,-1){15.3}}\put(20,11){$v$}

\put(4,18){$(X(A\oplus B))1\oplus (X(A)\oplus B))1$}
\put(6,17.8){\vector(1,-1){1.3}}\put(7,17){$(\mathcal{L}\otimes id)\oplus (\mathcal{L}\otimes id)$}

\put(17,20){$(X\oplus X)(A1)\oplus (X\oplus X)(B1)$}
\put(22,19.8){\vector(0,-1){17.3}}\put(22.2,13){$\mathcal{R}\oplus \mathcal{R}$}

\put(19.5,22){$(X\oplus X)(A1\oplus B1)$}
\put(22,21.8){\vector(0,-1){1.3}}\put(22.2,21){$\mathcal{L}$}
\put(22,22.5){\vector(0,1){3.3}}\put(22.2,24){$\mathcal{R}$}

\put(9,24){$X((A\oplus B)1)\oplus X((A\oplus B)1)$}
\put(14.5,24.5){\vector(2,1){2.7}}\put(10.5,25){$(id\otimes \mathcal{R})\oplus (id\otimes \mathcal{R})$}
\put(9,23.8){\vector(-1,-2){2.7}}\put(8,21){$a\oplus a$}

\put(0,26){$X((A\oplus B)1\oplus (A\oplus B)1)$}
\put(6,25.8){\vector(2,-1){2.8}}\put(7.9,25){$ \mathcal{L}$}
\put(6,26.5){\vector(2,1){2.8}}\put(7.9,27){$id\otimes( \mathcal{R}\oplus \mathcal{R})$}
\put(17.5,26){$X(A1\oplus B1)\oplus X(A1\oplus B1)$}
\put(18,25.8){\vector(-1,-3){2.45}}\put(14.5,21){$\mathcal{L}\oplus \mathcal{L}$}

\put(9,28){$X((A1\oplus B1)\oplus (A1\oplus B1))$}
\put(15,27.8){\vector(2,-1){2.8}}\put(17,27){$ \mathcal{L}$}

\put(12,26){(I)}
\put(4,23){(II)}
\put(12,21){(III)}
\put(7,14){(IV)}
\put(13.5,12){(V)}
\put(19,8){(VI)}
\put(9,9){(VII)}
\put(5,6){(VIII)}
\put(10,2){(IX)}
\put(18.5,21){(X)}
}
\put(11,-1.5){(2.2)}
\end{picture}
\end{center}

\begin{center}
\setlength{\unitlength}{0.53cm}
\begin{picture}(26,34)
\scriptsize{
\put(3.5,17){$X((A\oplus B)(1\oplus 1))$} 
\put(16,17){$(X\oplus X)(A1\oplus B1)$} 

\put(2.5,20){$X((A\oplus B)1\oplus(A\oplus B)1)$}
\put(15,20){$X(A1\oplus B1)\oplus X(A1\oplus B1)$}

\put(2.5,23){$X((A1\oplus B1)\oplus (A1\oplus B1))$}
\put(14,23){$(X(A1)\oplus X(B1))\oplus(X(A1)\oplus X(B1))$}

\put(2.5,26){$X(((A1\oplus B1)\oplus A1)\oplus B1)$}
\put(14,26){$((X(A1)\oplus X(B1)\oplus X(A1))\oplus X(B1)$}

\put(2.5,29){$X((A1\oplus (B1\oplus A1))\oplus B1)$}
\put(4,32){$X(A1\oplus (B1\oplus A1))\oplus X(B1)$}

\put(14,32){$(X(A1)\oplus X(B1\oplus A1))\oplus X(B1)$}

\put(14,29){$(X(A1)\oplus (X(B1)\oplus X(A1)))\oplus X(B1)$}

\put(2.5,14){$X(A(1\oplus 1)\oplus B(1\oplus 1))$}
\put(15,14){$(X\oplus X)(A1)\oplus (X\oplus X)(B1)$}

\put(14,11){$(X(A1)\oplus X(A1))\oplus (X(B1)\oplus X(B1))$}
\put(2.5,11){$X((A1\oplus A1)\oplus (B1\oplus B1))$}

\put(8,13){$X(A1\oplus A1)\oplus X(B1\oplus B1)$}
\put(2.5,8){$X(((A1\oplus A1)\oplus B1)\oplus B1)$}

\put(14,8){$((X(A1)\oplus X(A1))\oplus X(B1))\oplus X(B1)$}

\put(2.5,5){$X((A1\oplus (A1\oplus B1))\oplus B1)$}

\put(2.5,2){$X(A1\oplus (A1\oplus B1))\oplus X(B1)$}

\put(14,2){$(X(A1)\oplus X(A1\oplus B1))\oplus X(B1)$}

\put(14,5){$(X(A1)\oplus (X(A1)\oplus X(B1)))\oplus X(B1)$}}

\put(5.5,4.8){\vector(1,-1){2.3}}
\put(9.6,2.2){\vector(1,0){3.8}}
\put(16.7,2.7){\vector(1,1){2}}

\put(5.5,5.6){\vector(0,1){2}}%
\put(19,5.6){\vector(0,1){2}}

\put(5.5,10.6){\vector(0,-1){1.8}}
\put(19,10.6){\vector(0,-1){1.8}}
\put(5.5,13.6){\vector(0,-1){1.8}}
\put(19,13.6){\vector(0,-1){1.8}}
\put(5.5,16.6){\vector(0,-1){1.8}}
\put(19,16.6){\vector(0,-1){1.8}}
\put(5.5,17.8){\vector(0,1){1.8}}
\put(19,17.8){\vector(0,1){1.8}}
\put(5.5,20.8){\vector(0,1){1.8}}
\put(19,20.8){\vector(0,1){1.8}}
\put(5.5,23.8){\vector(0,1){1.8}}
\put(19,23.8){\vector(0,1){1.8}}

\put(5.5,28.7){\vector(0,-1){2}}
\put(19,28.7){\vector(0,-1){2}}

\put(5.5,29.7){\vector(1,1){2}}
\put(11,32.2){\vector(1,0){2.5}}
\put(17,31.8){\vector(1,-1){2}}

\put(7,11.5){\vector(1,1){1.3}}
\put(14,12.8){\vector(1,-1){1.3}}
\put(9,23){\vector(2,-1){5.5}}

\put(2.3,29){\line(-1,0){2.3}}
\put(0,29){\line(0,-1){24}}
\put(0,5){\vector(1,0){2.2}}
\put(23.2,29){\line(1,0){1.8}}
\put(25,29){\line(0,-1){24}}
\put(25,5){\vector(-1,0){1.8}}
\scriptsize{
\put(5.5,30.5){$\mathcal{L}$}

\put(11.7,32.5){$\mathcal{L}\oplus id$}%
\put(18.8,30.5){$(id\oplus \mathcal{L})\oplus id$}

\put(2.4,27.5){$id\otimes(a\oplus id)$}
\put(19.5,27.5){$a\oplus id$}

\put(4,24.5){$id\otimes a$}
\put(19.5,24.5){$a$}

\put(2.3,21.5){$id\otimes (\mathcal{R}\oplus \mathcal{R})$}
\put(19.5,21.5){$\mathcal{L}\oplus\mathcal{L}$}

\put(6,18.5){$id\otimes \mathcal{L}$}
\put(18,18.5){$\mathcal{R}$}

\put(3.5,15.5){$id\otimes \mathcal{R}$}
\put(19.5,15.5){$\mathcal{L}$}

\put(2.3,12.5){$id\otimes (\mathcal{L}\oplus \mathcal{L})$}
\put(19.5,12.5){$\mathcal{R}\oplus\mathcal{R}$}

\put(3.5,9.5){$id\otimes a$}
\put(19.5,9.5){$a$}

\put(2.3,6.5){$id\otimes(a\oplus id)$}
\put(19.5,6.5){$a\oplus id$}

\put(7,12){$\mathcal{L}$}
\put(15.3,12){$\mathcal{L}\oplus \mathcal{L}$}

\put(5.5,3.5){$\mathcal{L}$}
\put(11,1.5){$\mathcal{L}\oplus id$}%
\put(18.5,3.5){$(id\oplus \mathcal{L})\oplus id$}%

\put(0.2,18.5){$id\otimes((id\oplus c)\oplus id)$}
\put(21.5,18.5){$(id\oplus c)\oplus id$}
\put(12.5,21.7){$\mathcal{L}$}}

\put(11,26){(I)}
\put(11,17){(II)}
\put(11,8){(III)}
\put(1,17){(IV)}
\put(23,17){(V)}
\normalsize{\put(11,0){(2.3)}}
\end{picture}
\end{center}

\[\divide\dgARROWLENGTH by 2\begin{diagram} \node{X((A1\oplus(B1\oplus A1))\oplus
B1)}\arrow{s,lr}{\mathcal{L}}{\quad\quad\quad\quad\quad\quad\quad\quad\quad\text{
(I)}}\arrow[2]{e,t}{\quad id\otimes((id\oplus c)\oplus
id)\quad}\node[2]{X((A1\oplus(A1\oplus B1))\oplus
B1)}\arrow{s,r}{\mathcal{L}}\\
\node{X(A1\oplus(B1\oplus A1))\oplus X(
B1)}\arrow{s,lr}{id}{\quad\quad\quad\quad\quad\quad\quad\quad\quad\text{
(II)}}\arrow[2]{e,t}{\quad (id\otimes(id\oplus c))\oplus
id\quad}\node[2]{X(A1\oplus(A1\oplus B1))\oplus X(
B1)}\arrow{s,r}{\mathcal{L}\oplus id}\\
\node{(X(A1)\oplus X(B1\oplus A1))\oplus X(
B1)}\arrow{s,lr}{(id\oplus \mathcal{L})\oplus
id}{\quad\quad\quad\quad\quad\quad\quad\quad\quad\text{
(III)}}\arrow[2]{e,t}{\quad (id\oplus(id\otimes c))\oplus
id\quad}\node[2]{(X(A1)\oplus X(A1\oplus
B1))\oplus X(B1)}\arrow{s,r}{(id\oplus \mathcal{L})\oplus id}\\
\node{(X(A1)\oplus (X(B1)\oplus X(A1)))\oplus X(
B1)}\arrow[2]{e,t}{\quad (id\oplus c)\oplus
id\quad}\node[2]{(X(A1)\oplus (X(A1)\oplus X(B1)))\oplus X(B1)}
\end{diagram}\tag{2.4}\]
\end{proof}

\section{Categorical rings and the relationship with Ann-categories}
In [1], the authors presented the definition of categorical rings by modifying some axioms of the definition of Ann-categories. Let us recall this definition.

\begin{dn}
A categorical ring is a symmetric categorical group $\mathfrak{R}$ together with a bifunctor $\mathfrak{R}\times\mathfrak{R}\longrightarrow\mathfrak{R}$ (denoted by multiplication), an object $1\in\mathfrak{R}$, and natural isomorphisms:
\[a_{r, s, t} : (rs)t\longrightarrow r(st)\quad \text{(associative law),}\]
 \[\lambda_r: 1r\longrightarrow r ;\quad \rho_r: r1\longrightarrow
r\quad \text{(left and right unit),}\]
 \[\lambda_{r,x,y}:r(x+y)\longrightarrow rx+ry\quad \text{(left distributive law),}\]
 \[\rho_{x,y,s}:(x+y)s\longrightarrow xs+ys\quad \text{(right distributive law)}\]
such that $(\mathfrak{R},a, (1,\lambda,\rho))$ is a monoidal category making the diagrams  (1.1), (1.1'), (1.2),(1.3), (1.4), (1.4') and:
\[\begin{diagram}
\node[2]{r(x+y)+r(z+t)}\arrow{se,t}{\lambda_{r,x,y}+\lambda_{r,z,t}}\\
\node{r((x+y)+(z+t))}\arrow{ne,t}{\lambda_{r,x+y,z+t}}\arrow{s,l}{id\otimes
v_{x,y,z,t}}
\node[2]{(rx+ry)+(rz+rt)}\arrow{s,r}{v_{rx,ry,rz,rt}}\\
\node{r((x+z)+(y+t))}\arrow{se,b}{\lambda_{r,x+z,y+t}}\node[2]{(rx+rz)+(ry+rt)}\\
\node[2]{r(x+z)+r(y+t)}\arrow{ne,b}{\lambda_{r,x,z}+\lambda_{r,y,t}}
\end{diagram}\tag{3.1}\]
\[\begin{diagram}
\node[2]{(x+y)s+(z+t)s}\arrow{se,t}{\rho_{x,y,s}+\rho_{z,t,s}}\\
\node{((x+y)+(z+t))s}\arrow{ne,t}{\rho_{x+y,z+t,s}}\arrow{s,l}{v_{x,y,z,t}\otimes
id}
\node[2]{(xs+ys)+(zs+ts)}\arrow{s,r}{v_{xs,ys,zs,ts}}\\
\node{((x+z)+(y+t))s}\arrow{se,b}{\rho_{x+z,y+t,s}}\node[2]{(xs+zs)+(ys+ts)}\\
\node[2]{(x+z)s+(y+t)s}\arrow{ne,b}{\rho_{x,z,s}+\rho_{y,t,s}}
\end{diagram}\tag{3.1'}\]
commute.
\end{dn}
The main result of this paper is the relationship of Ann-categories and categorical rings. First, we have the following theorem.
\begin{thm}
Each Ann-category is a categorical ring.
\end{thm}

\begin{proof}
Assume that $(\A, \oplus, \otimes)$ is an Ann-category. We only need prove the commutation of the diagrams (3.1) and (3.1'). It can be deduced from the coherence theorem in an Ann-category [4]. However, we may present a direct proof as follows.\\
Consider the diagram:

\newpage 
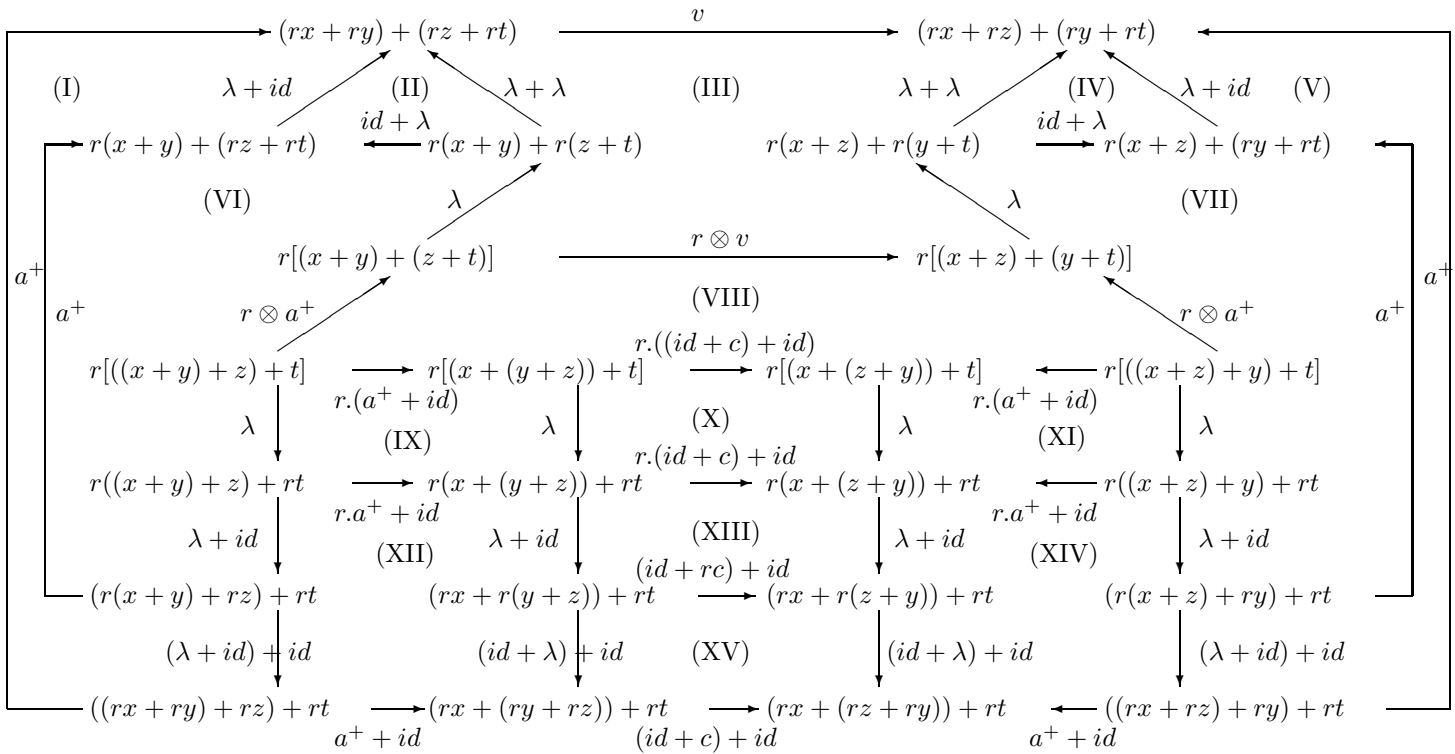
\begin{center}
\begin{sideways}
\setlength{\unitlength}{0.5cm}
\begin{picture}(32,24)
\put(0, 22){$(rx+ry)+(rz+rt)$}
\put(-5, 19){$r(x+y)+(rz+rt)$}
\put(0, 16){$r[(x+y)+(z+t)]$}
\put(-5, 13){$r[((x+y)+z)+t]$}
\put(-5, 10){$r((x+y)+z)+rt$}
\put(-5, 7){$(r(x+y)+rz)+rt$}
\put(-5, 4){$((rx+ry)+rz)+rt$}

\put(4, 19){$r(x+y)+r(z+t)$}
\put(4, 13){$r[(x+(y+z))+t]$}
\put(4, 10){$r(x+(y+z))+rt$}
\put(4, 7){$(rx+r(y+z))+rt$}
\put(4, 4){$(rx+(ry+rz))+rt$}

\put(13, 19){$r(x+z)+r(y+t)$}
\put(13, 13){$r[(x+(z+y))+t]$}
\put(13, 10){$r(x+(z+y))+rt$}
\put(13, 7){$(rx+r(z+y))+rt$}
\put(13, 4){$(rx+(rz+ry))+rt$}

\put(17, 22){$(rx+rz)+(ry+rt)$}
\put(22, 19){$r(x+z)+(ry+rt)$}
\put(17, 16){$r[(x+z)+(y+t)]$}
\put(22, 13){$r[((x+z)+y)+t]$}
\put(22, 10){$r((x+z)+y)+rt$}
\put(22, 7){$(r(x+z)+ry)+rt$}
\put(22, 4){$((rx+rz)+ry)+rt$}

\put(7.5, 22.2){\vector(1,0){9}}
\put(11, 22.5){$v$}

\put(25, 19.7){\vector(-3, 2){3}}
\put(24, 20.5){$\lambda+id$}
\put(18, 19.7){\vector(3, 2){3}}
\put(16.5, 20.5){$\lambda+\lambda$}

\put(7, 19.7){\vector(-3, 2){3}}
\put(6, 20.5){$\lambda+\lambda$}
\put(0, 19.7){\vector(3, 2){3}}
\put(-1.5, 20.5){$\lambda+id$}

\put(3.8, 19.2){\vector(-1,0){1.5}}
\put(2.2,19.6){$id+\lambda$}
\put(20.2, 19.2){\vector(1,0){1.5}}
\put(20.2,19.6){$id+\lambda$}

\put(4, 16.7){\vector(3, 2){3}}
\put(4.5, 17.5){$\lambda$}
\put(20, 16.7){\vector(-3, 2){3}}
\put(19.4, 17.5){$\lambda$}

\put(7.5, 16.2){\vector(1,0){9}}
\put(11, 16.5){$r\otimes v$}

\put(0, 13.7){\vector(3, 2){3}}
\put(-1, 14.5){$r\otimes a^+$}
\put(25, 13.7){\vector(-3, 2){3}}
\put(24, 14.5){$r\otimes a^+$}

\put(2, 13.2){\vector(1, 0){1.6}}
\put(1.5, 12.2){$r.(a^++id)$}
\put(11, 13.2){\vector(1, 0){1.6}}
\put(9.5, 13.7){$r.((id+c)+id)$}
\put(21.8, 13.2){\vector(-1, 0){1.6}}
\put(18.5, 12.2){$r.(a^++id)$}

\put(0,12.8){\vector(0,-1){2}}
\put(-1,11.5){$\lambda$}
\put(8,12.8){\vector(0,-1){2}}
\put(7,11.5){$\lambda$}
\put(16,12.8){\vector(0,-1){2}}
\put(16.5,11.5){$\lambda$}
\put(24,12.8){\vector(0,-1){2}}
\put(24.5,11.5){$\lambda$}

\put(2, 10.2){\vector(1, 0){1.6}}
\put(1.5, 9.2){$r.a^++id$}
\put(11, 10.2){\vector(1, 0){1.6}}
\put(9.5, 10.7){$r.(id+c)+id$}
\put(21.8, 10.2){\vector(-1, 0){1.6}}
\put(19, 9.2){$r.a^++id$}

\put(0,9.8){\vector(0,-1){2}}
\put(-2.4,8.5){$\lambda+id$}
\put(8,9.8){\vector(0,-1){2}}
\put(5.6,8.5){$\lambda+id$}
\put(16,9.8){\vector(0,-1){2}}
\put(16.4,8.5){$\lambda+id$}
\put(24,9.8){\vector(0,-1){2}}
\put(24.5,8.5){$\lambda+id$}

\put(11.2, 7.2){\vector(1, 0){1.5}}
\put(9.5, 7.7){$(id+rc)+id$}

\put(0,6.8){\vector(0,-1){2}}
\put(-3,5.5){$(\lambda+id)+id$}
\put(8,6.8){\vector(0,-1){2}}
\put(5.3,5.5){$(id+\lambda)+id$}
\put(16,6.8){\vector(0,-1){2}}
\put(16.2,5.5){$(id+\lambda)+id$}
\put(24,6.8){\vector(0,-1){2}}
\put(24.5,5.5){$(\lambda+id)+id$}

\put(2.5, 4.2){\vector(1, 0){1.4}}
\put(1.5, 3.2){$a^++id$}
\put(11.5, 4.2){\vector(1, 0){1.3}}
\put(9.5, 3.2){$(id+c)+id$}
\put(21.8, 4.2){\vector(-1, 0){1.2}}
\put(20, 3.2){$a^++id$}

\put(-5.2,7.2){\line(-1, 0){1}}
\put(-6.2,7.2){\line(0,1){12}}
\put(-6.2,19.2){\vector(1,0){1}}
\put(-5.9, 14.5){$a^+$}

\put(-5.2,4.2){\line(-1, 0){2}}
\put(-7.2,4.2){\line(0,1){18}}
\put(-7.2,22.2){\vector(1,0){7}}
\put(-7, 15.5){$a^+$}

\put(29.2,7.2){\line(1, 0){1}}
\put(30.2,7.2){\line(0,1){12}}
\put(30.2,19.2){\vector(-1,0){1}}
\put(29.2, 14.5){$a^+$}

\put(29.5,4.2){\line(1, 0){1.7}}
\put(31.2,4.2){\line(0,1){18}}
\put(31.2,22.2){\vector(-1,0){6.7}}
\put(30.5, 15.5){$a^+$}

\put(-6, 20.5){(I)}
\put(3, 20.5){(II)}
\put(11, 20.5){(III)}
\put(21, 20.5){(IV)}
\put(27, 20.5){(V)}

\put(-2, 17.5){(VI)}

\put(24, 17.5){(VII)}

\put(11, 14.9){(VIII)}

\put(2.8, 11.1){(IX)}
\put(11, 11.7){(X)}
\put(20.2, 11.2){(XI)}

\put(2.6, 8.1){(XII)}
\put(11, 8.7){(XIII)}
\put(20, 8.1){(XIV)}

\put(11, 5.5){(XV)}
\end{picture}
\end{sideways}
\end{center}
\newpage
In the above diagram, the regions (I) and (V) commute thanks to the naturality of $a^+$, the regions (II) and (IV) commute thanks to the composition of functors, the regions (VI), (VII), (XII) and (XIV) commute thanks to the compatibility of the functor $(L^a, \breve{L}^a)$ with the asociativity constraint $a^+$, the regions (IX), (X), (XI), (XIII) commute thanks to the naturality of the functor $\lambda$, the region (XV) commutes thanks to the compatibility of the functors $(L^a, \breve{L}^a)$ with the commutativity constraint $c$, the region (XIII) and the outside region commute thanks to the determination of the functor $v$ in the symmetric monoidal category $(\A, \oplus)$. Therefore, the region (III) commutes, that means the diagram (3.1) commutes.\\
With a similar proof, the diagram (3.1') commutes.

\end{proof}

\indent Now, to prove the converse, we need add the following axiom into the definition of categorical rings\\

(U) {\it For each object $a\in \mathfrak{R}$, the pairs $(L^a,\breve{L^a})$,and  $(R^a, \breve{R^a})$ defined by

\[
\begin{aligned}
&L^a & = &a\otimes- \;\;\;\; &R^a&=&-\otimes a\\
&\breve{L^a}_{x,y}&=&\mathcal{L}_{a, x, y}\;\;\;\; &\breve{R^a}_{x,
y}&=&\mathcal{R}_{x, y, a}
\end{aligned}
\]
are $\oplus$-functors which are compatible with the unitivity constraint $(0,g,d)$ of the operation $\oplus$. That means there exist isomorphisms $\widehat{L^A}:A\otimes O\longrightarrow O,\  \widehat{R^A}:O\otimes A\longrightarrow O, $ such that the diagrams (1.5), (1.5'), (1.6), (1.6') commute.\\

With this addition, we have the following theorem

\begin{thm}
Each categorical ring satisfying the condition (U) is an Ann-category.
\end{thm}

\begin{proof}

Assume that $\A$ is a categorical ring satisfying the condition (U). We must show that $\A$ satisfies the axiom (Ann-1) of an Ann-category. According to Proposition 2, it remains to show that the functor $(L^a, \breve{L}^a)$ is compatible with the associativity constraint $a^+$, i.e, the commutation of the following diagrams:

\[\begin{diagram}
\node{x(a+(b+c))}\arrow{e,t}{\lambda}\arrow{s,r}{id\otimes
a}\node{xa+x(b+c)}\arrow{e,t}{id\oplus \lambda}\node{xa+(xb+xc)}\arrow{s,r}{a}\\
\node{x(a+b)+c)}\arrow{e,t}{\lambda}\node{x(a+b)+xc}\arrow{e,t}{\lambda\oplus id}
\node{(xa+xb)+xc}
\end{diagram}\tag{3.2}\]
and a similar diagram for the pair $(R^A,\breve{R}^A),$ for each A.

First, since $(L, \breve{L^A})$ is compatible with the constraint $(0, g, d)$, there exists a functor $\widehat{L}^A,$ such that these diagrams:

\[\begin{diagram}
\node{A(0+X)}\arrow{e}\arrow{s}\node{A0+AX}\arrow{s}\\
\node{AX}\arrow{e}\node{0+AX}
\end{diagram}\quad\quad
\begin{diagram}
\node{A(X+0)}\arrow{e}\arrow{s}\node{AX+A0}\arrow{s}\\
\node{AX}\arrow{e}\node{AX+0}
\end{diagram}
\]
commute.\\
In order to prove the diagram (3.2), let's consider the diagram

\begin{center}
\setlength{\unitlength}{0.5cm}
\begin{picture}(26,20)
\footnotesize{
\put(0,0){$(XA+XC)+(X0+XD)$}
\put(18.5,0){$(XA+X0)+(XC+XD)$}
\put(5,3){$(XA+XC)+(0+XD)$}
\put(15,3){$(XA+0)+(XC+XD)$}
\put(5,6){$(XA+XC)+XD$}
\put(15,6){$XA+(XC+XD)$}
\put(0,9){$X(A+C)+X(0+D)$}
\put(19.5,9){$X(A+0)+X(C+D)$}
\put(5,12){$X(A+C)+XD$}
\put(16,12){$XA+X(C+D)$}
\put(5,15){$X((A+C)+D)$}
\put(16,15){$X(A+(C+D))$}
\put(0,18){$X((A+C)+(0+D))$}
\put(19.5,18){$X((A+0)+(C+D))$}

\put(7.3,0.2){\vector(1,0){10.7}}
\put(5,0.8){\vector(1,1){2}}
\put(22,0.8){\vector(-1,1){2}}
\put(11.7,3.2){\vector(1,0){2.8}}
\put(8,5.8){\vector(0,-1){2}}
\put(19,5.8){\vector(0,-1){2}}
\put(1,8.8){\vector(0,-1){8}}
\put(25,8.8){\vector(0,-1){8}}
\put(5,8.8){\vector(1,-1){2}}
\put(22,8.8){\vector(-1,-1){2}}
\put(8,11.8){\vector(0,-1){5}}
\put(19,11.8){\vector(0,-1){5}}
\put(1,17.8){\vector(0,-1){8}}
\put(25,17.8){\vector(0,-1){8}}
\put(8,14.8){\vector(0,-1){2}}
\put(19,14.8){\vector(0,-1){2}}
\put(5,9.8){\vector(1,1){2}}
\put(22,9.8){\vector(-1,1){2}}
\put(5,17.8){\vector(1,-1){2}}
\put(22,17.8){\vector(-1,-1){2}}
\put(6.1,18.2){\vector(1,0){13}}
\put(10.3,6.2){\vector(1,0){4.3}}
\put(9.6,15.2){\vector(1,0){6}}}
\scriptsize{
\put(13,0.4){$v$}
\put(2,1.5){$id+(\widehat{L}+id)$}
\put(21.5,1.5){$(id+\widehat{L})+id$}
\put(13,3.4){$v$}
\put(1.2,4.5){$\rho+\rho$}
\put(6.1,4.5){$id+g$}
\put(19.2,4.5){$d+id$}
\put(23.3,4.5){$\rho+\rho$}
\put(13,6.4){$a^+$}
\put(3,7.5){$\rho+id\otimes g$}
\put(21.2,7.5){$id\otimes d+\rho$}
\put(8.2,9){$\rho+id$}
\put(17,9){$id+\rho$}
\put(2.5,10.5){$id+id\otimes g$}
\put(21.5,10.5){$id\otimes d+id$}
\put(1.2,13.5){$\rho$}
\put(24.5,13.5){$\rho$}
\put(7.5,13.5){$\rho$}
\put(18.5,13.5){$\rho$}
\put(12,15.4){$id\otimes a^+$}
\put(2.5,16.5){$id\otimes(id+g)$}
\put(21.2,16.5){$id\otimes(d+id)$}
\put(12,18.4){$id\otimes v$}
}
\put(4,13.5){(I)}
\put(6,10.2){(II)}
\put(4,4.5){(III)}
\put(12.5,16.5){(IV)}
\put(12.5,10.5){(V)}
\put(12.5,4.5){(VI)}
\put(12.5,1.5){(VII)}
\put(21.5,13.5){(VIII)}
\put(19.5,10.2){(IX)}
\put(22,6){(X)}
\end{picture}
\end{center}

In the above diagram, the outside region commutes thanks to the hypothesis (1.5); the regions (I) and (VIII) commute thanks to the functorical property of $\rho$, the regions (II) and (IX) commute thanks to the composition, the regions (III) and (X) commute thanks to the definition of $\widehat{L}$; the regions (IV) and (VI) commute thanks to the coherence theorem in a symmetric monoidal category, the region (VII) commutes thanks to the 
functorical property of $v$. Therefore, the region (V) commutes. In other words, the diagram (3.2) commutes.

The compatibility of the functor $(R^A,\breve{R}^A)$ with the associativity constraint $a^+$ can be proved similarly.

\end{proof}

\begin{center}

\end{center}
\end{document}